\documentclass{article}

\usepackage{amssymb,amsmath,eufrak,amsthm,amscd}

\newcommand{\ord}{\operatorname{ord}}
\newcommand{\degtr}{\operatorname{tr\,deg}}
\newcommand{\diffdegtr}{\operatorname{tr\,deg^\Delta}}

\newcommand{\diffspec}{\operatorname{Spec^\Delta}}

\newcommand{\Hom}{\operatorname{Hom}}

\newcommand{\Qt}{\operatorname{Qt}}

\newcommand{\height}{\operatorname{ht}}
\newcommand{\diffheight}{\operatorname{ht_\Delta}}
\newcommand{\diffdim}{\operatorname{dim_\Delta}}
\newcommand{\type}{\operatorname{type}}
\newcommand{\rk}{\operatorname{rk}}

\newtheorem{theorem}{Theorem}

\newtheorem{corollary}[theorem]{Corollary}
\newtheorem{statement}[theorem]{Statement}
\newtheorem*{statement*}{Statement}
\newtheorem*{theorem*}{Theorem}
\newtheorem*{lemma*}{Lemma}
\newtheorem*{fact*}{Fact}

\theoremstyle{definition}

\newtheorem*{definition*}{Definition}
\newtheorem{example}[theorem]{Example}

\newtheorem*{example*}{Example}

\theoremstyle{remark}

\newtheorem*{remark*}{Remark}

\author{D.\,V.~Trushin}
\title{Local dimension of differential algebraic variety}
\date{}

\begin{document}

\maketitle

\begin{abstract}
We consider a relation between local and global characteristics of a
differential algebraic variety. We prove that dimension of tangent
space for every regular point of an irreducible differential
algebraic variety coincides with dimension of the variety.
Additionally, we classify tangent spaces at regular points in the
case of one derivation.
\end{abstract}

\section{Introduction}

We investigate the local behavior of differential algebraic variety.
We mostly attend to tangent spaces and relation between their
characteristics and characteristics of the variety. We prove two
main results. The first one gives an estimate on dimension of
associated graded algebra (theorem~\ref{maintheorem}). The second
result is a classification of tangent spaces at regular points in
the case of one derivation (theorem~\ref{tang}).

In section~\ref{sec2} we briefly recall some definitions. Next
section~\ref{sec3} is entirely devoted to an auxiliary technique
used further. In subsection~\ref{sec31} the notion of height of
finitely generated algebra over a field is discussed. In
subsection~\ref{seec32} we introduce the notion of differential
dimension. Differential dimension is defined for differentially
finitely generated algebras over a field (section~\ref{sec321}) and
for differential algebras of finite type (section~\ref{sec322}). To
investigate the local structure of differential algebraic variety we
need the notion of tangent space. In section~\ref{sec33} we
introduce the notion of linear differential space. All tangent
spaces appeared are linear differential spaces. So, investigation of
such spaces is a preparation for investigation of tangent spaces.
Subsection~\ref{sec331} is devoted to the definition and basic
properties. Additionally, a theorem about correspondence between
submodules and subspaces is proven in this subsection
(theorem~\ref{Nullth}). Subsection~\ref{sec332} introduces some
useful technique for further applications. Using this machinery, in
subsection~\ref{sec333}, we define differential dimension polynomial
of linear differential space. Section~\ref{sec4} is devoted to the
first main result -- theorem~\ref{maintheorem}. This result gives us
an estimate on dimension of associated graded algebra. Additionally,
we present an example showing that in this estimate the inequality
can not be changed to the equality. Next section~\ref{sec5} is
devoted to an application of obtained results to tangent spaces. In
subsection~\ref{sec52} we recall the notion of regular point. We
show that differential dimension polynomial of tangent space at
regular point coincides with differential dimension polynomial of
the variety (statement~\ref{dimst}). The result of
section~\ref{sec4} is applied in subsection~\ref{sec53} to describe
the behavior of points similar to regular ones. In last
section~\ref{sec6} we classify tangent spaces at regular points of
an irreducible differential algebraic variety in the case of one
derivation. The main result is theorem~\ref{tang}. Additionally, we
present an example showing that structure of a tangent space depends
on the choice of a regular point.

\paragraph{Acknowledgement.} The author appreciate to Marina
Kondrateva for her remarks and advices. The correct version of
theorem~\ref{tang} and example~\ref{examp} appeared due to Marina
Kondrateva.

\section{Terms and notation}\label{sec2}

The word ring means an associative commutative ring with an identity
element. All homomorphisms preserve the identity element. A
differential ring is a ring with a finitely many pairwise commuting
derivations. The set of all derivations will be denoted by $\Delta$.
Additionally, we suppose differential rings to be Ritt algebras. The
field of fractions of an integral domain $A$ will be denoted by
$\Qt(A)$.

\section{Preliminaries}\label{sec3}

In this section we develop an auxiliary technique.

\subsection{Dimension}\label{sec31}

The section is devoted to the notion of the height of finitely
generated algebra over a field. Additionally, we shall fix notation
and recall the notion of dimension.

Let $K$ be a field, and let $B$ be a finitely generated algebra over
$K$. Let a family $x=(x_1,\ldots,x_n)$ generate the algebra $B$ over
$K$. Define the following subspaces
$$
I={<}x_1,\ldots,x_n{>}_K,\:\: I^0=K
$$
and
$$
I_n=\sum_{k=0}^n I^k.
$$
Consider the following function $\chi_x(n)=\dim_K I_n$.

\begin{statement}
Using notation above, the following holds:
\begin{enumerate}
\item The function $\chi_x(n)$ coincides with a polynomial for sufficiently large $n$.
\item If $y=(y_1,\ldots,y_k)$ is another system of generators, then
$$
\deg\chi_x=\deg\chi_y.
$$
\end{enumerate}
\end{statement}
\begin{proof}

Let us show the first statement. The algebra $B$ can be presented as
follows
$$
K[x_1,\ldots,x_n]/\frak a.
$$
Let some order on monomials of $K[x_1,\ldots,x_n]$ be fixed. We
suppose that the order preserve degree, for example deglex. Let
$g_1,\ldots,g_s$ be a Gr\"obner basis of $\frak a$. Then dimension
of $I_n$ coincides with the number of all monomials of degree less
then of equal to $n$ such that they are not divided by the leading
monomials of $g_i$. Due to lemma~\cite[chapter~0, sec.~17,
lemma~16]{Kl} it is clear that the number of such monomials
coincides with some polynomial.

Let us prove the second statement. Let some other system of
generators $y$ be given. The corresponding sequence of subspaces
will be denoted by $J_n$. Then for some $m_0$ we have the inclusion
$J\subseteq I_{m_0}$. Hence
$$
J^m\subseteq(I_{m_0})^m=\left(\sum_{k=0}^{m_0}I^k\right)^m\subseteq\sum_{k=0}^{m_0
m}I^k=I_{m_0 m}.
$$
Consequently, $\chi_{y}(m)\leqslant\chi_x(m_0 m)$. In analogue way
we have that for some $n_0$ the inequality
$\chi_{x}(n)\leqslant\chi_y(n_0 n)$ holds. Since both functions are
polynomials, $\deg\chi_x=\deg\chi_y$.
\end{proof}

Degree of polynomial $\chi_x$ will be called a height of algebra $B$
and will be denoted by $\height B=\deg \chi_x$. The previous
statement guaranty that the notion of height does not depend on the
set of generators. Now we shall prove some basic properties of
height.

\begin{statement}
Let $A$ and $B$ be finitely generated algebras over a field $K$.
Then the following holds:
\begin{enumerate}
\item If $A\subseteq B$, then $\height A\leqslant \height B$.
\item If $B$ is a quotient algebra of algebra $A$, then  $\height A \geqslant \height
B$.
\item If $A\subseteq B$ and $B$ is integral over $A$, then  $\height A=\height
B$.
\item The height of polynomial ring in $n$ variables coincides with $n$, in other words
$$
\height K[x_1,\ldots,x_n]=n.
$$
\item Let some set of algebraically independent elements $y_1,\ldots,y_n$ of algebra $A$ be
chosen such that algebra $A$ is integral over $K[y_1,\ldots,y_n]$,
then  $\height A=n$.
\end{enumerate}
\end{statement}
\begin{proof}
(1). Let $x=(x_1,\ldots,x_n)$ be a family of generators of $A$. This
family can be extended to a set of generators of $B$. Assume that
$x'=x\cup\{\,x_{n+1},\ldots,x_{r}\,\}$. Let $I={<}x_1,\ldots,x_n{>}$
and $J={<}x_1,\ldots,x_r{>}$. Then $I_n\subseteq J_n$. Consequently,
$\chi_x(n)\leqslant\chi_{x'}(n)$. And thus $\deg \chi_x\leqslant\deg
\chi_{x'}$.

(2). Let $x=(x_1,\ldots,x_n)$ be a family of generators of $A$. The
images of elements of $x$ will be chosen as the generators of $B$.
Let $I_n$ and $J_n$ be corresponding sequences of subspaces in $A$
and $B$ respectively. Then $J_n$ is the image of $I_n$ under the
quotient mapping. Consequently,
$\chi_x(n)\geqslant\chi_{\overline{x}}(n)$. So, $\deg
\chi_x\geqslant\deg\chi_{\overline{x}}$.

(3). From item~(1) it follows that $\height A\leqslant \height B$.
We just need to show the other inclusion. Let a family
$m_1,\ldots,m_r$ generate $B$ as a module over $A$. Then for some
$x_{ij}^k$ we have $m_im_j=\sum_k x_{ij}^km_k$. Let us chose a
family of generators of $A$ such that the set $x_{ij}^k$ is included
to the chosen family of generators. Let denote this family by $x$
and the corresponding sequence of subspaces by $I_n$. We shall
denote by $y$ the following set of generators of $B$: $x_im_j$
whenever $x_i\in x$ $1\leqslant j\leqslant r$. The corresponding
sequence of subspaces for $y$ will be denoted by $J_n$. From the
definition we have $J=Im_1+\ldots+Im_r$. Then
\begin{gather*}
J_n=\sum_{k=0}^{n}J^k=\sum_{k=0}^{n}\left(\sum_{i=1}^rIm_i\right)^k=\sum_{k=0}^{n}\sum
I^km_{i_1}\ldots m_{i_k}\subseteq\\
\subseteq\sum_{k=0}^{n}\sum_{i=1}^rI^{2k-1}m_i=\sum_{i=1}^r\left(\sum_{k=0}^{n}I^{2k-1}\right)m_i\subseteq
\sum_{i=1}^r I_{2n-1}m_i.
\end{gather*}
Consequently, $\chi_{y}(n)\leqslant r\chi_{x}(2n-1)$. So,
$\deg\chi_y\leqslant\deg\chi_x$.

(4). A straightforward calculation show that the number of all
monomials in $n$ variables of degree not greater than $m$ is equal
to $\binom{n+m}{n}$. The last expression is a polynomial in $n$
variables of degree $m$.

(5). The statement is an immediate corollary of previous two items.
\end{proof}

Now we shall recall the notion of dimension for a ring $A$. We use
the same notation as in books~\cite{AM} and~\cite{Mu}. The Krull
dimension of algebra $A$ will be denoted by $\dim A$. If $A$ is a
local ring with a maximal ideal $\frak m$, then $G_{\frak m} (A)$
will denote the corresponding associated graded ring, in other words
$$
G_{\frak m} (A)=\bigoplus_{k=0}^\infty \frak m^n/\frak m^{n+1}.
$$
Let
$$
P(t)=\sum_{n=0}^\infty \dim_{A/\frak m} \frak m^n/\frak m^{n+1}
$$
be the Poincar\'e series of $G_{\frak m} (A)$ (see~\cite[chapter~11,
sec.~1]{AM}). Order of the pole of $P(t)$ at $t=1$ we shall denote
by $d(A)$. Speaking about height of the algebra $G_{\frak
m}(A_{\frak m})$, we shall consider this ring as the algebra over
the field $A/\frak m$.

For every local noetherian ring $A$ there is the equality $\dim
A=d(A)$ (see~\cite[chapter~11, sec.~2, th.~11.14]{AM}). From the
other hand, if $A$ is a finitely generated integral domain over a
field, then it is universally catenary (see~\cite[chapter~5,
sec.~14, coroll.~3(2)]{Mu}). Thus for every maximal ideal $\frak m$
of $A$ we have $\dim A=d(A_{\frak m})$.

\begin{statement}\label{heightstate}
Let $A$ be a finitely generated algebra over a field $K$, and let
$\frak m$ be a maximal ideal of $A$. Then the following holds:
\begin{enumerate}
\item $\height A=\dim A$.
\item If $A$ is an integral domain, then
$$
\height G_{\frak m}(A_{\frak m})=d(A_{\frak m})=\dim A.
$$
\end{enumerate}
\end{statement}
\begin{proof}
(1). Consider the Noether's normalization of $A$. So, there is a
family of algebraically independent elements $y_1,\ldots,y_r$ of $A$
such that $A$ is integral over $K[y_1,\ldots,y_r]$. Then from
condition~(5) of previous statement it follows that $\height A=r$.
Corollary~\cite[chapter~5, sec.~14, coroll.~1]{Mu} guaranties that
$\dim A=r$.

(2). Corollary~\cite[chapter~5, sec.~14, coroll.~3(1)]{Mu} shows
that $d(A_{\frak m})=\dim A$. Corollary~\cite[chapter~11, sec.~1,
coroll.~11.2]{AM} guaranties the other equality.
\end{proof}

\subsection{Differential dimension}\label{seec32}

\subsubsection{Differentially finitely generated algebras}\label{sec321}

At present there are two different ways to define differential
dimension of differentially finitely generated algebra over a
differential field. We shall recall both definitions and show that
they gives the same results.

Let $K$ be a differential field with a set of derivations $\Delta$,
and let the number of derivations equals $m$. Consider an algebra
$A$ differentially finitely generated over a field $K$.  Let
$x=(x_1,\ldots,x_n)$  be a system of differential generators of $A$,
in other words $A=K\{x_1,\ldots,x_n\}$. We shall define the notion
of differential height of $A$. Let $A_k=K[\theta_1
x_1,\ldots,\theta_n x_n\mid \ord \theta_i\leqslant k]$. For every
ideal $\frak a$ of $A$ the image of the system $x$ in $A/\frak a$
will be denoted by $\overline{x}$. Consider the function
$\chi^A_{x}(k)=\dim A_k$.

\begin{statement}\label{diffheight}
Using notation above the following holds:
\begin{enumerate}
\item The function $\chi^A_{x}(t)$ is a polynomial for sufficiently large $t$.

\item Let $l$ be degree of the polynomial $\chi^A_x(t)$ and $a_k$ be its coefficient
at $t^k$. Then the numbers  $l$ and $a_l$ do not depend on the
system $x$. Moreover, $l$ does not exceed $m$, and the number
$d_l=l!a_l$ is integer. Particulary, the number $d_m=m! a_m$ is
integer.

\item Let $\frak p_1,\ldots, \frak p_r$ be the set of all minimal prime
(and thus differential) ideals of $A$. Then
$$
\chi^A_{x}(k)=\max_{1\leqslant i\leqslant r} \chi^{A/\frak
p_i}_{\overline{x}}(k)=\chi^{A/\frak p_j}_{\overline{x}}(k)
$$
for some $j$.

\item If $A$ is an integral domain, then $d_m$ coincides with
differential transcendence degree of the fraction field of $A$ over
$K$.
\end{enumerate}
\end{statement}
\begin{proof}
Let $\frak n$ be the nilradical of $A$, then its contraction to
$A_k$ coincides with the nilradical $\frak n_k$. Since $\dim
A_k=\dim A_k/\frak n_k$, we may suppose that our algebra  is
reduced. Let contraction of ideal $\frak p_i$ to $A_k$ be denoted by
$\frak p_{ik}$. Then from the definition of Krull dimension we have
$\dim A_k=\max_{i}\dim A_k/\frak p_{ik}$. Since $A/\frak p_i$  is an
integral domain, theorem~\cite[chapter~2, sec.~12, th.~6(a)]{Kl}
guaranties that for sufficiently large $k$ the functions
$$
\chi^{A/\frak p_i}_{\overline{x}}(k)=\dim(A/\frak p_i)_k=\dim
A_k/\frak p_{ik}
$$
are polynomials of degree not greater than $m$. Consequently,
$$
\max_{1\leqslant i\leqslant r} \chi^{A/\frak p_i}_{\overline{x}}(k)
$$
coincides with one of the polynomials $\chi^{A/\frak
p_j}_{\overline{x}}(k)$ for some $j$. The last polynomial coincides
with $\chi_{x}(k)$. So, we have proven~(1) and~(3).

Let us show that degree $l$ and the leading coefficient $a_l$ do not
depend on the choice of $x$. Let $y=(y_1,\ldots,y_s)$  be  other
system of differential generators. The corresponding filtration will
be denoted by $A'_k$. Then for some $k_0$ we see that $x_i$ belong
to $A'_{k_0}$. Consequently, we have $A_{k}\subseteq A'_{k+k_0}$.
And thus
$$
\chi^A_x(k)\leqslant\chi^A_y(k+k_0).
$$
By the similar arguments we have the other inequality
$$
\chi^A_y(r)\leqslant\chi^A_x(r+r_0).
$$
So, degree and leading coefficient of these polynomials coincide.
The fact that $l$ is not greater than $m$ has been proven.  Let us
show that $l!a_l$ is integer. Indeed, since  $\chi^A_x(t)$ is
integer for sufficiently large $t$, the polynomial
$$
\Delta\chi^A_x(t)=\chi^A_x(t)-\chi^A_x(t-1).
$$
has the same property. Then the number
$d_l=\Delta^l\chi^A_x(t)=l!a_l$ is integer. The last item of the
statement follows from  theorem~\cite[chapter~2, sec.~12,
th.~6(c)]{Kl}.
\end{proof}

We shall use the notation of the previous statement. The number $l$
is called a type of differential algebra $A$ and denoted by $\type
A$. The number $d_l$ is called a typical differential height of $A$.
The number $d_m$ is called a differential height of $A$ and is
denoted by $\diffheight A$. This definition is due to J.~Kovacic. If
$A$ is a local differential ring of finitely generated type over a
field ($A$ is a localization of differentially finitely generated
algebra), saying about differential height of $G_{\frak m}(A)$, we
shall consider this ring as an algebra over the field $A/\frak m$.

The following definition appeared in~\cite[sec.~1,
pp.~207-208]{Jh01}. For every differentially finitely generated
algebra over a field $K$ consider the function
$$
\mu\colon \diffspec A\times\diffspec A\to \mathbb Z\cup\{\infty\}
$$
defined in lemma~\cite[pp.~208]{Jh01}. Then the maximal value of the
function $\mu$ will be called a differential Krull type of $A$ and
denoted by $\type_K A$. Theorem~\cite[sec.~2, pp.~208]{Jh01} says
that differential Krull type does not exceed $m$. Consider a longest
chain of prime differential ideals $\frak p_0\supseteq \frak
p_1\supseteq \ldots\supseteq \frak p_k$ such that $\mu(\frak
p_i,\frak p_{i+1})=\type_K A$. Then the maximum of such $k$ will be
called a typical Krull dimension of $A$. Differential Krull
dimension of $A$ is the maximum of $k$ such that $\frak p_0\supseteq
\frak p_1\supseteq \ldots\supseteq \frak p_k$ and $\mu(\frak
p_i,\frak p_{i+1})=m$. It is clear, that differential Krull
dimension equals zero if and only if differential Krull type is less
than $m$. Differential Krull dimension will be denoted by $\diffdim
A$. It should be noted that our definition of Krull dimension a
little bit differ from the initial one. Theorem~\cite[sec.~2,
pp.~208]{Jh01} says that if $A$ is an integral domain, then
differential Krull dimension coincides with differential
transcendence degree of fraction field of $A$ over $K$. From
item~(4) of the previous statement we have.

\begin{statement}
Let $A$ be a differentially finitely generated integral domain over
a field $K$. Then $\diffheight A=\diffdim A$.
\end{statement}

We shall show that the equality holds in the case of arbitrary
differentially finitely generated algebra.

\begin{statement}\label{diffhtdim}
For any algebra $A$ differentially finitely generated over $K$ there
is the equality $\diffheight A=\diffdim A$.
\end{statement}
\begin{proof}
From statement~\ref{diffheight} it follows that $\diffheight A$
coincides with the maximum of $\diffheight A/\frak p_i$, where
$\frak p_i$ is a minimal prime differential ideal in $A$. From the
other hand, from definition it follows that $\diffdim A$ coincides
with the maximum of $\diffdim A/\frak p_i$, where $\frak p_i$ is a
minimal prime differential ideal of $A$. From previous statement it
follows that both these numbers coincide.
\end{proof}

\subsubsection{Differential algebras of finite type}\label{sec322}

Let $A$ be a differential algebra over a field $K$ such that there
exists a system of elements $x=(x_1,\ldots,x_n)$ in $A$ with the
following property. Algebra $A$ coincides with
$K\{x_1,\ldots,x_n\}_\frak p$, where $\frak p$ is a prime ideal of
$K\{x_1,\ldots,x_n\}$. The system of elements $x$ will be called a
system of local generators of $A$. Let $B$ denote the differential
subalgebra in $A$ generated by $x_1,\ldots,x_n$. The maximal ideal
of $A$ will be denoted by $\frak m$. The images of $x_i$ in $B/\frak
p$ will be denoted by $\overline{x}_i$. Then the system
$\overline{x}=(\overline{x}_1,\ldots,\overline{x}_n)$ is a system of
differential generators of $B/\frak p$. For a given system of
generators $x$ we shall define the following rings
$B_r=K[\theta_ix_i\mid \ord\theta_i\leqslant r]$ and $(B/\frak
p)_r=K[\theta_i\overline{x}_i\mid \ord\theta_i\leqslant r]$. Then
define  $\frak p_r=B_r\cap\frak p$, $A_r=(B_r)_{\frak p_r}$, and the
maximal ideal of $A_r$ will be denoted by $\frak m_r$. Consider the
function $\chi^A_x(t)=\dim A_t$.

\begin{statement}
Using notation above, the following holds:
\begin{enumerate}
\item The function $\chi^A_{x}(t)$ coincides with a polynomial for sufficiently large
$t$.

\item Let $l$ be degree of the polynomial $\chi^A_x(t)$, and $a_k$ be its coefficients at
$t^k$. Then numbers $l$ and $a_l$ do not depend on choice of $x$.
Moreover,  $l$ is not greater than  $m$, and the number $d_l=l!a_l$
is integer. Particulary, the number $d_m=m! a_m$ is integer.

\item Let $\frak p_1,\ldots, \frak p_r$ be all minimal prime (and thus differential)
ideals of algebra $A$. Then
$$
\chi^A_{x}(k)=\max_{1\leqslant i\leqslant r} \chi^{A/\frak
p_i}_{\overline{x}}(k)=\chi^{A/\frak p_j}_{\overline{x}}(k)
$$
for some $j$.
\end{enumerate}
\end{statement}
\begin{proof}

Let us show the equality $\dim A_t=\dim B_t-\dim (B/p)_t$. Since $B$
can be embedded to $A$, every minimal prime ideal of $A$ is
contained in ideal $\frak p$. Let us denote the contraction of
$\frak p_i$ to $B_t$ by $\frak p_{it}$. Then from
corollary~\cite[chapter~14, sec~14.H, coroll.~3]{Mu} we have
\begin{gather*}
\dim A_t=\max_{1\leqslant i\leqslant r}( \dim A_t/\frak
p_{it})=\max_{1\leqslant i\leqslant r}(\dim B_t/\frak p_{it})-\dim B_t/\frak p_t=\\
=\dim B_t-\dim B_t/\frak p_t=\dim B_t-\dim (B/\frak p)_t.
\end{gather*}
So,
$$
\chi^A_x(t)=\max_{1\leqslant i\leqslant r}\chi^{A/\frak
p_i}_{\overline{x}}(t)=\chi^B_{x}(t)-\chi^{B/\frak
p}_{\overline{x}}(t).
$$
Additionally, the polynomial $\chi^A_x(k)$ has integer values and is
of degree not greater than $m$, because it is a difference of two
polynomials of degree not greater than $m$.

Let $y=(y_1,\ldots,y_s)$ be some other system of local generators of
$A$. The corresponding filtration will be denoted by $A'_k$. Then
for some $k_0$ we have $A_0\subseteq A'_{k_0}$, and thus
$A_{k}\subseteq A'_{k+k_0}$. Hence
$$
\chi^A_x(k)\leqslant\chi^A_y(k+k_0).
$$
Using analogous arguments we have the other inequality
$$
\chi^A_y(t)\leqslant\chi^A_x(t+t_0).
$$
Therefore degrees and the leading coefficients of these polynomials
coincide.

Since polynomial $\chi_x(k)$ has integer values for sufficiently
large $k$, the polynomial $\Delta \chi_x(t)=\chi_x(t)-\chi_x(t-1)$
has integral values too. So, the coefficient
$d_l=\Delta^l\chi_x(t)=l!a_l$  is integer, where $a_l$ is the
leading coefficient.
\end{proof}

Degree of the polynomial $\chi^A_x$ is called a differential type of
algebra $A$. Coefficient $d_m$ is called a differential height of
$A$, and coefficient $d_l$ is called a typical differential height.

\subsection{Linear differential spaces}\label{sec33}

In this section we prepare a basic technique for local theory.

\subsubsection{Definition}\label{sec331}

Let $K$ be a differential field. Concider an affine space $K^n$. A
differential module $K[\Delta]^n$ will be denoted by $\mathcal L_n$.
For every element $\xi\in\mathcal L_n$ we define the mapping
$\xi\colon K^n\to K$ by the following rule. Let
$\xi=(\xi_1,\ldots,\xi_n)$, where  $\xi_i\in K[\Delta]$, and let
$x=(x_1,\ldots,x_n)\in K^n$. Then $\xi(x)=\xi_1x_1+\ldots+\xi_nx_n$.
The elements of module $\mathcal L_n$ will be called linear
differential functions on $K^n$ and will be denoted by $\mathcal
L(K^n)$.

Consider an arbitrary differential submodule $N\subseteq\mathcal
L_n$. For this submodule we define the following set
$$
V(N)=\{\,x\in K^n\mid \forall \xi\in N\colon \xi(x)=0  \,\}.
$$
Conversely, for every subset $X\subseteq K^n$ we define a
differential submodule
$$
I(X)=\{\,\xi\in \mathcal L_n\mid \xi|_X=0\,\}.
$$
The sets of the form $V(N)$ will be called linear differential
spaces.

For every natural number $n$ consider the set of differential
homomorphisms $\Hom_{K[\Delta]}(K[\Delta]^n,K)$. Then there is the
mapping
$$
\Hom_{K[\Delta]}(K[\Delta]^n,K)\to K^n
$$
by the rule $\xi\mapsto (\xi(e_1),\ldots,\xi(e_n))$, where  $e_i$ is
a standard basis of $K[\Delta]^n$. Since $K[\Delta]^n$ is free
module, the constructed mapping is bijective. Let now $N\subseteq
K[\Delta]^n$ be a submodule. Then the set of all homomorphisms of
$\Hom_{K[\Delta]}(K[\Delta]^n,K)$ vanishing on $N$ maps to $V(N)$
bijectively. It is clear that these homomorphisms can be identified
with $\Hom_{K[\Delta]}(K[\Delta]^n/N,K)$. In other words we have the
bijection
$$
\Hom_{K[\Delta]}(K[\Delta]^n/N,K)\to V(N)
$$
by the rule $\xi\mapsto (\xi(e_1),\ldots,\xi(e_n))$.

Now consider the symmetric algebra on the module $K[\Delta]^n$. So,
$$
R_n=S_K(K[\Delta]^n).
$$
If $e_1,\ldots,e_n$ are standard free generators of the module
$K[\Delta]^n$, then the ring $R_n$ coincides with a differential
polynomial ring $K\{e_1,\ldots,e_n\}$. Then we may suppose that
module $K[\Delta]^n$ is embedded into $R_n$. For every submodule
$N\subseteq K[\Delta]^n$ we define the ideal $[N]=(N)$. Since $[N]$
is a graded ideal with respect to degree, we have $[N]\cap
K[\Delta]^n=N$. The ideal $[N]$ is prime because it is generated by
linear differential polynomials. Every differential homomorphism
$\xi\in\Hom_{K[\Delta]}(K[\Delta]^n,K)$ gives a differential
homomorphism $\xi\colon R_n\to K$. The last correspondence is a
bijection between the set of all differential homomorphisms of
$\Hom_{K[\Delta]}(K[\Delta]^n,K)$ and the set of all differential
rings homomorphisms $R_n\to K$.

\begin{statement}[Nullstellensatz]\label{Nullth}
For every differential submodule $N\subseteq \mathcal L_n$ we have
$N=I(V(N))$.
\end{statement}
\begin{proof}

The inclusion $N\subseteq I(V(N))$ is obvious. Let us show the other
one. Consider the algebra $R_n$ and the ideal $[N]$. Let $x\in
K[\Delta]^n\setminus N$. Since $[N]\cap K[\Delta]^n=N$, $x\notin
[N]$. We know that the ideal $[N]$ is prime and the field $K$ is
differentially closed. Therefore there exists a point $a\in K^n$
such that for every  $f\in [N]$ we have $f(a)=0$ and $x(a)\neq0$.
So, $x$ is not in $I(V(N))$.
\end{proof}

Let $V\subseteq K^n$ be a linear differential space. Consider the
restriction of functions in $K[\Delta]^n$ to $V$. The resulting
module will be denoted by $\mathcal L(V)$ and called the module of
linear differential functions on $V$. From the previous theorem it
follows that $\mathcal L(V)$ is isomorphic to $K[\Delta]^n/I(V)$.
Then the symmetric algebra $S_K(\mathcal L(V))$ can be identified
with coordinate ring of differential algebraic variety $V\subseteq
K^n$.

Let $\varphi_1,\ldots,\varphi_m$ be elements of $K[\Delta]^n$. Then
we define a linear differential mapping $\varphi\colon K^n\to K^m$
by the rule: if $x\in K^n$, the coordinates of the point
$\varphi(x)$ equal $\varphi_1(x),\ldots,\varphi_m(x)$. Let $V$ and
$U$ be linear differential spaces in $K^n$ and $K^m$ respectively.
Then the mapping  $\varphi\colon V\to U$ is called a linear
differential mapping if it coincides with a restriction of some
linear differential mapping from $K^n$ to $K^m$.

\subsubsection{Characteristic sets}\label{sec332}

Let $W=K[\Delta]^n$ be a free differential module, and let its
standard basis be denoted by $e_1,\ldots,e_n$. Then the module $W$
is generated by $\theta e_i$, where $\theta\in \Theta$. By a ranking
of $e_1,\ldots,e_n$ we shall mean a total ordering of the set of all
derivatives $\theta e_i$ that satisfies the two conditions
$$
u\leqslant\theta u,\:\: u\leqslant v\Rightarrow \theta u\leqslant
\theta v.
$$
From statement~\cite[chapter~0, sec.~17, lemma~15]{Kl} it follows
that ranking exists and every ranking is a well ordering of the set
of all derivatives $\theta e_i$. A ranking will be said to be
orderly if the rank of $\theta e_i$ is less then that of $\theta'
e_j$ whenever $\ord\theta<\ord\theta'$.

Let $w$ be an element of $W$. Then it is a linear combination of
elements $\theta e_i$. Let $u_w$ denotes a leading derivative
involved in $w$. So, there appear a pre-order on the set $W$.
Namely, we compare two elements by their leading vectors. Let $w'$
be other element of $W$. We shall say that $w'$ is reduced with
respect to $w$ if $w'$ is free of every derivative of $u_w$.  The
set of elements of $W$ will be called reduced if every element of
this set is reduced with respect to every other element. The fact
that every autoreduced set is finite is a corollary of
lemma~\cite[chapter~0, sec.~17, lemma~15(a)]{Kl} is


Now we shall define a particular order on the set of all autoreduced
sets of $W$. Let $A=\{f_1,\ldots,f_r\}$ and $B=\{g_1,\ldots,g_s\}$
be autoreduced sets such that their elements are arranged in order
of increasing rank. We shall say that rank of $A$ is less than rank
of $B$ if
\begin{enumerate}
\item There exist $k\in \mathbb N$ $k\leqslant r$ and $k\leqslant
s$ such that  $u_{f_i}=u_{g_i}$ whenever $1\leqslant i<k$ and
$u_{f_k}<u_{g_k}$.
\item Or $r>s$ and $u_{f_i}=u_{g_i}$ whenever $1\leqslant i\leqslant s$.
\end{enumerate}
It should be noted that if $r=s$ and for every $i$ we have
$u_{f_i}=u_{g_i}$ $1\leqslant i\leqslant s$, then $A$ and $B$ has
the same rank.

\begin{statement}
In every nonempty set of autoreduced subsets of $W$ there exists an
autoreduced set of minimal rank.
\end{statement}
\begin{proof}

Let $\mathcal A$ be any nonempty set of autoreduced subsets of $W$.
Define by induction an infinite decreasing sequence of subsets of
$\mathcal A$ by the conditions that $\mathcal A=\mathcal A_0$ and,
for $i>0$, $\mathcal A_i$ is the set of all autoreduced sets $A\in
\mathcal A_{i-1}$ with $|A|\geqslant i$ such that the $i$-th lowest
element of $A$ is of lowest possible rank. It is obvious that in all
elements of $\mathcal A_i$ the $i$-th lowest elements have the same
leader $v_i$. If every $\mathcal A_i$ were nonempty, then the
leaders $v_i$ would form an infinite sequence of derivatives of
$e_i$ such that no $v_i$ is a derivative of any other, and this
would contradict lemma~\cite[chapter~0, sec.~17, lemma~15(a)]{Kl}.
Therefore there is a smallest $i$ such that $\mathcal A_i=\emptyset$
and, since $\mathcal A=\mathcal A_0\neq\emptyset$, $i>0$. Any
element of $\mathcal A_{i-1}$ is clearly an autoreduced subset in
$\mathcal A$ of lowest rank.
\end{proof}

For every submodule $N\subseteq W$ there exists a minimal
autoreduced subset consisting of elements of $N$. Such autoreduced
sets we shall call characteristic sets of $N$.

\subsubsection{Dimension}\label{sec333}

There are two equivalent methods to define dimension of linear
differential space. Let $V$ be a linear differential space, and let
$\mathcal L(V)$ be the set of all linear differential functions on
$V$. Using the module $\mathcal L(V)$ we shall define dimension of
$V$.

The set of elements  $y_1,\ldots, y_d$ of differential module $M$
over a field $K$ is called differentially independent if the set
$\theta y_i$ is linearly independent. The maximal number of
differentially independent elements of module $M$ is called a
differential dimension of $M$.

The second way of defining differential dimension is the following.
Let $y_1,\ldots,y_n$ be differential generators of $M$. The there is
a family of subspaces in $M$
$$
M_k=<\theta_1 y_1,\ldots,\theta_n y_n\mid \ord \theta_i\leqslant
k>_K.
$$
All subspaces $M_k$ are of finite dimension. Therefore there is a
function $\varphi(k)=\dim_K M_k$. Statement~\cite[chapter~I, sec.~4,
th.]{Jh02} guaranties that for sufficiently large $k$ the function
$\varphi(k)$ coincides with a polynomial of degree not greater that
$|\Delta|$.  This polynomial is called a differential dimension
polynomial. The coefficient at term $t^m$ has the form $a_m/m!$.
Then the number $a_m$ is called differential dimension of $M$.
Statement~\cite[chapter~III, sec~2, prop.]{Jh02} says that these
both definitions coincide. We shall denote differential dimension of
module $M$ by  $\diffdim_K M$. Moreover
statement~\cite[chapter~III,sec.~2, lemma]{Jh02} says that dimension
has an appropriate behavior. Namely, for every exact sequence of
differential nodules $0\to M'\to M\to M''\to 0$, there is the
equality $\diffdim_K M=\diffdim_K M'+\diffdim_K M''$. Differential
dimension of linear differential space is a differential dimension
of its module of linear differential functions.

The next our purpose is to derive the information about other
coefficients of differential dimension polynomial in the case of one
derivation.  In the following statement we suppose that we deal with
the case of one derivation $\delta$.

\begin{statement}\label{freeterm}
Let $M$ be a differentially finitely generated module over a field
$K$ of differential dimension $d$. Let $m_1,\ldots,m_n$ be a system
of differential generators, and let $\varphi(t)=dt+r$ be the
differential dimension polynomial calculated by the given set of
generators. Then among the elements
$$
m_1,\ldots,m_n
$$
there are $d$ differentially independent elements
$m_{i_1},\ldots,m_{i_d}$ such that the quotient module
$M/[m_{i_1},\ldots,m_{i_d}]$ is a vector space of dimension $r$ over
$K$.
\end{statement}
\begin{proof}

The module $M$ can be presented as follows $M=K[\delta]^n/N$, where
standard basis $e_i$ maps to $m_i$. Let us fix an orderly ranking.
And let $F$ be a characteristic set for $N$. Changing the order of
elements $m_i$ we can suppose that $m_1,\ldots,m_d$ are
differentially independent and the leaders of elements of $F$ are
derivatives of $m_{d+1},\ldots,m_n$. Assume that leaders have the
following form  $\theta_{d+1} e_{d+1},\ldots,\theta_n e_n$. It is
clear that in this situation a free term of differential dimension
polynomial coincides with number of elements in the following set
$$
\{\,\theta e_i \mid d+1\leqslant i\leqslant n,\: \theta e_i<
\theta_i e_i \,\}.
$$
But this set form a basis of the quotient space
$$
M/[m_1,\ldots,m_d].
$$
\end{proof}

\section{Associated graded algebra}\label{sec4}

Let $A$ be a differential algebra over a field $K$ such that there
exists a family of elements $x=(x_1,\ldots,x_n)$ with the following
property. The algebra $A$ coincides with $K\{x_1,\ldots,x_n\}_\frak
p$, where $\frak p$ is a prime differential ideal of
$K\{x_1,\ldots,x_n\}$. Such system of elements $x$ we shall call a
system of local generators of $A$. Let $B$ denote the algebra
generated over $K$ by elements $x_1,\ldots,x_n$. The maximal ideal
of $A$ will be denoted by $\frak m$. The images of $x_i$ in $B/\frak
p$ will be denoted by $\overline{x}_i$. Then the system
$\overline{x}=(\overline{x}_1,\ldots,\overline{x}_n)$ is the family
of differential generators of $B/\frak p$. For our system of
generators we define the rings $B_r=K[\theta_ix_i\mid
\ord\theta_i\leqslant r]$ and $(B/\frak
p)_r=K[\theta_i\overline{x}_i\mid \ord\theta_i\leqslant r]$. Then we
define  $\frak p_r=B_r\cap\frak p$, $A_r=(B_r)_{\frak p_r}$, and the
maximal ideal of $A_r$ will be denoted by $\frak m_r$.

Consider an associated graded algebra
$$
G_{\frak m}(A)=A/\frak m\oplus \frak m/\frak
m^2\oplus\ldots\oplus\frak m^n/\frak m^{n+1}\oplus\ldots
$$
It coincides with a direct limit of associated graded algebras
$G_{\frak m_r}(A_r)$. Moreover, $G_{\frak m}(A)$ is differentially
finitely generated over $A/\frak m$.

\begin{theorem}\label{maintheorem}
Let notation be as above. If $B$ is an integral domain, then
\begin{enumerate}
\item $\type G_{\frak m}(A)\leqslant \type A$.

\item If $\type G_{\frak m}(A)= \type A$, then typical differential height
of  $G_{\frak m}(A)$ is not greater than typical differential height
of $A$.

\item $\diffdim G_{\frak m}(A)\leqslant\diffdim B-\diffdim B/\frak p=\diffheight A$.
\end{enumerate}
\end{theorem}
\begin{proof}
Corollary~\cite[chapter~5, sec.~14, coroll.~3(1)]{Mu} says that
$$
\dim B_r=\dim A_r+\dim(B_r/\frak p_r).
$$
We know that $B_r/\frak p_r=(B/\frak p)_r$.
Theorem~\cite[chapter~11, sec.~2, th.~11.14]{AM} with
statement~\ref{heightstate} guaranty that there is the following
sequence of equalities
$$
\dim A_r=d(A_r)=d(G_{\frak m_r}(A_r))=\height (G_{\frak m_r}(A_r)).
$$
Consider the natural mapping  $G_{\frak m_r}(A_r)\to G_{\frak
m}(A)$. Let us denote its image by $D'_r$. Consider the smallest
$A/\frak m$ algebra generated by $D'_r$. In other word, consider
$A/\frak m\cdot D'_r$. We shall denote this algebra by $D_r$. If
$$
D'_r=\oplus_{k\geqslant0} D'_{rk},
$$
where $D'_{r0}=A_r/\frak m_r$, then
$$
D_r=\oplus_{k\geqslant0} D_{rk},
$$
where $D_{rk}=A/\frak m\cdot D'_{rk}$. Consequently, dimension of
$D'_{rk}$ over $K$ coincides with dimension of $D_{rk}$ over
$A/\frak m$. So,
$$
\height D_r=\height D'_r\leqslant \height G_{\frak m_r}(A_r).
$$
From the construction we have  $\delta(D_r)\subseteq D_{r+1}$ for
every $r$ and every $\delta\in \Delta$. Moreover,
$$
G_{\frak m}(A)=\bigcup_{k\geqslant 0} D_k.
$$
Let $y=(y_1,\ldots,y_t)$ be a system of differential generators of
$G_{\frak m}(A)$ over $A/\frak m$. Then for some $r_0$ it follows
that for all $i$ we have $y_i\in D_{r_0}$. Therefore
$$
G_{\frak m}(A)_r\subseteq D_{r+r_0}.
$$
Hence we have the following sequence of inequalities
\begin{gather*}
\dim G_{\frak m}(A)_r\leqslant \dim D_{r+r_0}=\height
D_{r+r_0}\leqslant\\
\leqslant \height G_{\frak m_{r+r_0}}(A_{r+r_0})=\dim B_{r+r_0}-\dim
(B/\frak p)_{r+r_0}.
\end{gather*}
So, for some $r_0$ we have that for all $r$ the following holds
$$
\dim G_{\frak m}(A)_r\leqslant\dim B_{r+r_0}-\dim (B/\frak
p)_{r+r_0}=\dim A_{r+r_0}.
$$
Since the functions above are polynomials with integer values, all
three results immediately follows from the last inequality and
definitions.
\end{proof}

We shall show that in the previous theorem we can not change the
inequality to an equality. The following example appeared in
work~\cite{Jh04}. We just add some remarks to it.

\begin{example}
Let $K$ be a field with one differential operator and  $K{<}y{>}$ be
a differential extension of $K$ such that $y$ and $y'$ are
algebraically independent over $K$ and $y''=y'/y$. If we let
$B=K\{y\}$, then $B_0=K[y]$ and $B_r=K[y,y'/y^{r-1}]$ if $r>0$.
Clearly $By$ is a prime differential ideal of $B$, and in fact
$B/By=K$. Let $A=B_{By}$. Since $y'=y^{r-1}(y'/y^{r-1})\in By^{r-1}$
for all $r>1$, $y'\in\cap_{d\geqslant 0}\frak m^d$. Let us show that
$[y']=\cap_{d\geqslant 0}\frak m^d$. Indeed, It is easy to see that
the ring $A'=A_{By}/[y']$ coincides with  $K[y]_{(y)}$, where
$y'=0$. But the ring $A'$ is regular, so the intersection of powers
of maximal ideal of $A'$ is zero. Let us note that under the
homomorphism $A\to A'$ the powers of maximal ideal of $A$
corresponds to the powers of maximal ideal of $A'$ (because the
kernel is contained in the intersection $\cap_{d\geqslant 0}\frak
m^d$). It is clear that  $G_{\frak m}(A)=G_{\frak m}(A')=K[t]$. So,
typical differential height of $A$ is equal to its transcendence
degree $2$, but typical differential height of associated graded
ring is $1$.
\end{example}

\section{Cotangent spaces}\label{sec5}

\subsection{Definition}\label{sec51}

Let $K$ be a differential field with the set of derivations
$\Delta$. Let $B$ be a differentially finitely generated algebra
over $K$ not containing nilpotent elements. Then for such algebra
$B$ we can produce a differential algebraic variety as follows. Let
$y_1,\ldots,y_n$ be a family of differential generators of $B$. Then
$B=K\{y_1,\ldots,y_n\}/\frak a$, where $\frak a$ is a radical
differential ideal. We define the set $X=V(\frak a)$ consisting of
all common zeros for $\frak a$ in $K^n$. So, $B$ can be identified
with a coordinate ring of $X$ is sense of~\cite{Cs}. Let $\frak
m\subseteq B$ be a maximal differential ideal of $B$.  Since $K$ is
differentially closed, $\frak m=[y_1-a_1,\ldots,y_n-a_n]$. So all
maximal differential ideals correspond to points of $X$. If $x\in X$
is a point of the variety, then the corresponding maximal
differential ideal will be denoted by $\frak m_x$ or simpler by
$\frak m$, when it is known what point is under consideration.

We shall define cotangent space for a given point $x$. Consider
local ring $A=B_\frak m$, where $\frak m$ corresponds to $x$. Then
the ideal $\frak m$ is differentially finitely generated in $A$.
Therefore the quotient module $\frak m/\frak m^2$ is differentially
finitely generated over $A/\frak m$. Since $K$ is differentially
closed the field $A/\frak m$ coincides with $K$. In other words,
module $\frak m/\frak m^2$ is differentially finitely generated over
$K$. This module will be called a cotangent space for $x$ and
denoted by $T^*_x$.

We shall describe the second construction of cotangent space. Let
$\Omega_{A/K}$ be the module of K\"ahler differentials of $A$ over
$K$. In~\cite{Jh03} it is shown that this module is a differential
module over $A$. Then statement~\cite[chapter~10, sec.~25,
th.~58]{Mu} guaranties that $\Omega_{A/K}\mathop{\otimes}_{A}A/\frak
m=\frak m/\frak m^2$. It is easy to see that the last isomorphism is
an isomorphism of differential modules. This definition allows to
calculate cotangent spaces in applications. Namely, if the
differential ring $B$ is of the following form
$$
K\{y_1,\ldots,y_n\}/[f_1,\ldots,f_s],
$$
then module of differentials has the following form
$$
\Omega_{B/K}=<dy_1,\ldots,dy_n>/[df_1,\ldots,df_s].
$$
Using localization and tensor products, we are able to calculate
cotangent space of every point.

Tangent space for $x$ is a linear differential space corresponding
to module $\frak m/\frak m^2$. Explicitly, let $y_1,\ldots,y_n$  be
the set of differential generators of $\frak m/\frak m^2$. Then this
module has the following form $\frak m/\frak m^2 =K[\Delta]^n/N$.
The tangent space $T_x$ coincides with $V(N)\subseteq K^n$.

\subsection{Regular points of differential spectrum}\label{sec52}

We shall use the definition of a regular point of differential
spectrum from~\cite{Jh04}. Whenever a local ring $A$ is given we
shall always associate with it the following notation
\begin{itemize}
\item $\frak m$ the maximal ideal of $A$.
\item $P(A)=S_K(\frak m/\frak m^2)$ the symmetric algebra over
$K$ on vector space $\frak m/\frak m^2$.
\item $G(A)=\oplus_{d\geqslant 0}\frak m^d/\frak m^{d+1}$ the associated graded algebra of $A$.
\item $\tau_A\colon P(A)\to G(A)$ the unique $K$-algebra homomorphism that extends the identity mapping of
$\frak m/\frak m^2$.
\end{itemize}

Let $A$ be a differential algebra over $K$ such that there is a
system of elements $x=(x_1,\ldots,x_n)$ in $A$ with the following
property. The algebra $A$ coincides with $K\{x_1,\ldots,x_n\}_\frak
p$, where $\frak p$ is a prime differential ideal of
$K\{x_1,\ldots,x_n\}$. The $K$-algebra generated by elements
$x_1,\ldots,x_n$ will be denoted by $B$. The maximal ideal of $A$
will be denoted by $\frak m$. For our system of generators we define
the ring $B_r=K[\theta_ix_i\mid \ord\theta_i\leqslant r]$. Then we
define $\frak p_r=B_r\cap\frak p$, $A_r=(B_r)_{\frak p_r}$, and the
maximal ideal of $A_r$ will be denoted by $\frak m_r$.

We shall say that the ring $A$ is regular if the following
conditions hold:
\begin{description}
\item{({\bf A1}):} The mapping $\tau_A\colon P(A)\to G(A)$ is an isomorphism.
\item{({\bf A2}):} There exists a system of local generators for $A$
such that
$$
\phi_r\colon A/\frak m\mathop{\otimes}_{A_r}\Omega_{A_r/K}\to
A/\frak m\mathop{\otimes}_{A}\Omega_{A/K}
$$
is injective for all $r\geqslant 0$.
\end{description}

Let $B$ be a differentially finitely generated algebra over a field
$K$. Let $\frak p$ be a prime differential ideal of $B$. We shall
say that $\frak p$ is a regular point of differential spectrum if
the local ring $B_\frak p$ is  regular. The set of local generators
in condition ({\bf A2}) we shall call the set of regular local
generators of $A$ or the set of regular generators for $\frak p$ in
$B$.

It should be noted that for every differentially finitely generated
integral domain the set of all regular points is open and everywhere
dense~\cite[sec.~5, pp.~228, th.]{Jh04}.

For every system of local generators $x=(x_1,\ldots,x_n)$ of $A$ the
set $dx=(dx_1,\ldots,dx_n)$ is a system of differential generators
of $\Omega_{A/K}$. Therefore their images
$\overline{dx}=(\overline{dx}_1,\ldots,\overline{dx}_n)$  are
differential generators of
$$
A/\frak m\mathop{\otimes}_A\Omega_{A/K}=\frak m/\frak m^2 .
$$
Define the sequence of modules $M_k\subseteq M=\frak m/\frak m^2$ as
follows
$$
M_k=<\theta_i \overline{dx}_i\mid \ord\theta_i \leqslant k>_K.
$$
Let $\Omega'_r=A\mathop{\otimes}_{A_r}\Omega_{A_r/K}$. Then it is
easy to see that
$$
M_r=\phi_r(A/\frak m\mathop{\otimes}_A\Omega'_r).
$$
Dimension of vector space $M_t$ over $K$ will be denoted by
$\chi^{M}_{\overline{dx}}(t)$.

\begin{statement}\label{dimst}
Let notation be as above. Suppose that the local ring $A$ is
regular. Let $B\subseteq A$ be a differential subalgebra generated
by regular system of generators. The ideal $\frak p$ will denote the
contraction of $\frak m$ to $B$. Then
$$
\chi^{M}_{\overline{dx}}(t)=\chi^B_x(t)=\chi^{B/\frak
p}_{\overline{x}}(t)+\chi^A_x(t).
$$
\end{statement}
\begin{proof}

Since $A$ is regular, then $A$ is an integral domain~\cite[sec.~2,
prop.~1]{Jh04}. Since finitely generated integral domain $B_t$ is
universally catenary, then we have the equality
$\chi^B_x(t)=\chi^{B/\frak p}_{\overline{x}}(t)+\chi^A_x(t)$
(see~\cite[chapter~14, sec.~14.H, coroll.~3]{Mu}).

Since $A_t$ is a regular local ring (see~\cite[sec.~2,
prop.~1]{Jh04}), then
$$
\degtr_K \Qt(A)=\rk \Omega_{A_t/K}
$$
(see~\cite[sec.~1, lemma~4(2)]{Jh04}). Consequently, we have the
sequence of equalities
$$
\dim B_t=\degtr_K A=\rk \Omega_{A_t/K}=\dim_K (A/\frak
m\mathop{\otimes}_{A_t}\Omega_{A_t/K}).
$$
Since the mapping $\phi_r$  is an isomorphism onto its image, then
$$
\dim_K (A/\frak m\mathop{\otimes}_{A_t}\Omega_{A_t/K})=\dim_K M_t.
$$
We have $\dim B_t=\dim_K M_k$. So,
$\chi^B_x(t)=\chi^M_{\overline{dx}}(t)$.
\end{proof}

\begin{corollary}\label{dimcor}

Let $B$ be a differentially finitely generated integral domain over
a field $K$. And let $\frak p$ be a regular point of differential
spectrum of $B$. Define $A=B_\frak p$ and $\frak m$ is its maximal
ideal. Then
$$
\diffdim_K \frak m/\frak m^2=\diffdim B-\diffdim B/\frak p.
$$
\end{corollary}

\subsection{Dimension}\label{sec53}

Let $B$ be a differentially finitely generated algebra over a field
$K$, and let $\frak p$ be a prime differential ideal of $B$.
Consider the ring $A_\frak p$. Its maximal ideal will be denoted by
$\frak m$. We are interested in local rings satisfying condition
({\bf A1}) from the definition of regular point of differential
spectrum.

\begin{theorem}\label{geomain}
Let $B$ be a differentially finitely generated integral domain over
$K$ and let $\frak p$ be a point of differential spectrum satisfying
condition \emph{({\bf A1})}. Let us denote the residue field of
$\frak p$ by $L$, the fraction field of $B$ by $F$, and $A=B_\frak
p$ be a local ring with a maximal ideal $\frak m$. Then the
following holds
$$
\diffdim \frak m/\frak m^2=\diffdim B -\diffdim B/\frak p
$$
or using the language of the fields
$$
\diffdim\frak m/\frak m^2=\diffdegtr_{K}F -\diffdegtr_{K}L
$$
\end{theorem}
\begin{proof}
The inequality
$$
\diffdim \frak m/\frak m^2\geqslant\diffdim B -\diffdim B/\frak p
$$
was proven in theorem~\cite[sec.~4, pp.~96]{Jh03}. We shall show the
other one. Indeed, for every given point $\frak m$ we know that
$$
P(A)=G(A).
$$
Then, from one hand, differential dimension of module $\frak m/\frak
m^2$ coincides with differential dimension of algebra
$P(A)=S_K(\frak m/\frak m^2)$. From the other hand, from
theorem~\ref{maintheorem} it follows that differential dimension of
$G(A)$ is not greater than
$$
\diffdim B -\diffdim B/\frak p,
$$
Q.E.D.
\end{proof}

The geometrical meaning of the previous theorem is the following.
For every point of differential algebraic variety satisfying
condition ({\bf A1}) differential dimension of tangent space
coincides with differential dimension of the variety. If we combine
this result with theorem~\cite[sec.~5, pp.~228]{Jh04} we obtain that
the set of all such points contains an open everywhere dense set.

\section{Classification}\label{sec6}

The last attainment of our work is a classification of tangent
spaces at regular points of differential algebraic varieties in the
case of one derivation. The detailed theory of differential
algebraic varieties is given in~\cite{Cs}. We shall borrow terms and
notation from this work.

Let $K$ be a differentially closed field with the subfield of
constants $C$. Let $X$ be an irreducible differential algebraic
variety over $K$. Let $B$ be a coordinate ring for $X$.

\begin{theorem}\label{tang}
Let $X$ be an irreducible differential algebraic variety. Let
$x=(x_1,\ldots,x_n)$ be coordinates on $X$. So,
$$
B=K\{X\}=K\{x_1,\ldots,x_n\}/I(X).
$$
Assume that a point $y\in X$ is regular in coordinates
$x_1,\ldots,x_n$, and let differential dimension polynomial has the
following form $\chi^B_x(t)=dt+r$. Then the tangent space at $y$ is
of the form $T_y=K^d\times C^k$, where $d$ coincides with
differential dimension of $X$ and $k$ is less then or equal to $r$.
\end{theorem}
\begin{proof}

Let $\frak m$ be the maximal ideal of a point $y\in X$. Consider the
cotangent space $T^*_y=\frak m/\frak m^2$. This module is a left
module over the ring of differential operators $K[\delta]$. The last
ring is left and right euclidian~\cite[chapter~2, sec.~1,
lemma~2.1]{SW}. Since $\frak m/\frak m^2$ is finitely generated
module over $K[\delta]$, it is isomorphic to
$$
K[\delta]^{\oplus n}\oplus K[\delta]/p_1K[\delta]\oplus\ldots\oplus
K[\delta]/p_sK[\delta].
$$
We shall denote the submodule
$K[\delta]/p_1K[\delta]\oplus\ldots\oplus K[\delta]/p_sK[\delta]$ by
$V$. Then $V$ is a vector space of finite dimension over $K$. We
shall denote its dimension by $k$. Let $e_1,\ldots,e_k$ be a basis
of $V$ over $K$. Then
$$
\delta(e_1,\ldots,e_k)=(e_1,\ldots,e_k)C
$$
For some matrix $C\in M_k(K)$. Let us change the basis using a
matrix $B$ as follows
$$
(e'_1,\ldots,e'_k)=(e_1,\ldots,e_k)B.
$$
So,
$$
\delta(e'_1,\ldots,e'_k)=(e'_1,\ldots,e'_k)C'
$$
Then the matrix $C'$ coincides with $B^{-1}CB+B^{-1}\delta B$. We
shall solve the equation $C'=0$. The equation is equivalent to
$\delta B=-CB$. Statement~\cite[chapter~1, sec.~3, prop.~1.20]{SW}
says that there exists a desired nondegenerated matrix $B$ with
coefficients in some Picard-Vessiot extension. But the field $K$ is
differentially closed, therefore it contains all Picard-Vessiot
extensions. So, the matrix $B$ belongs to $GL_k(K)$. Using vectors
$e'_1,\ldots,e'_k$ as a basis, we have
$$
V=(K[\delta]/\frak \delta K[\delta])^{\oplus k}.
$$
So, the module $\frak m/\frak m^2$ is isomorphic to
$$
K[\delta]^{\oplus n}\oplus(K[\delta]/\frak \delta K[\delta])^{\oplus
k}.
$$
By the definition it is clear that $n=\diffdim \frak m/\frak m^2$.
And from corollary~\ref{dimcor} it follows that $n=\diffdim B=\dim
X$. The corresponding tangent space has the form $K^d\times C^k$.

We shall show the estimate for $k$. Statement~\ref{dimst} says that
differential dimension polynomial of $\frak m/\frak m^2$ calculated
in induced coordinates coincides with $dt+r$. Submodule $V$
intersects any free submodule $T$ of $T^*_y$ by zero. From
statement~\ref{freeterm} it follows that there is a free submodule
$T$ in $T^*_y$ such that $T^*_y/T$ has dimension $r$. But $V$ can be
embedded to $T^*_y/T$, and therefore $k=\dim_K V\leqslant\dim_K
T^*_y/T=r$
\end{proof}

\begin{example}\label{examp}
The goal of this example is to show that coefficient $k$ from the
previous statement depends on the choice of a regular point and can
take any value from $0$ to a value of free term of differential
dimension polynomial.

Let $K$ be a differential field with a subfield of constants $C$.
Consider the ring of differential polynomials $K\{z,y\}$ and the
ideal $[zy'-y]$. Let us note that the ring $K\{z,y\}_z/[zy'-y]$ has
no zero divisors. From criterion~\cite[sec.~3, pp.~219, th.]{Jh04}
it follows that points  $(z,y)$ such that $z\neq 0$ are regular
points of the variety $X$ given by the equation $zy'-y=0$. Let us
note that there is a unique point with condition $z=0$. This point
is $(0,0)$.

Consider tangent spaces in all points except $(0,0)$. We shall show
that for every point $(z_0,y_0)$ such that $z_0\neq 0, y'_0\neq 0$
the tangent space is $K$, and for other points (conditions are
$z_0\neq 0, y'_0=0$) the tangent space is $K\times C$. Indeed,
cotangent space at point $\frak m=[z-z_0,y-y_0]$ has the following
form
$$
\frak m/\frak m^2=\Omega_{A/K}\otimes_{A}A/\frak m,
$$
where  $A$ is local ring corresponding to the point. Then
$$
\frak m/\frak m^2={<}dz,dy{>}/[z_0dy'+y'_0dz-dy].
$$
The last module can be presented as
$$
K[\delta]\oplus K[\delta]/(y'_0,-1+z_0\delta).
$$
Note if $y'_0=0$, then cotangent space is $K[\delta]\oplus
K[\delta]/\delta K[\delta]$ otherwise $K[\delta]$.

We shall show the examples of these two possibilities. For the point
$(1,0)$ tangent space is specified by the equation $dy'-dy=0$. Let
$\gamma$ be a nonzero solution of the last equation. Then the set of
all points $(z,c\gamma)$, where $z\in K$ and $c\in C$, is the
tangent space at $(1,0)$ and is isomorphic to $K\times C$. Consider
the point $(t,t)$, where $t\in K$ such that $t'=1$. Then the tangent
space is specified by the equation $tdy'+dz-dy=0$. Hence the set of
all points $(y-ty',y)$, where $y\in K$, is the tangent space at
point $(t,t)$ and is isomorphic to $K$.
\end{example}


\begin{thebibliography}{99}

\bibitem{AM} M.\,F.~Atiyah, I.\,G.~Macdonald.  Introduction to commutative
algebra.  Addison-Wesley. 1969.

\bibitem{Cs} P.\,J.\,Cassidy. Differential algebraic groups. Amer.
J.~Math., 1972, 93, n.~3, pp.~891-954.

\bibitem{Jh01} J.\,Johnson. A notion of Krull dimension for differential
rings. Comment. Math. Helv., 44, 1969, 207-216.

\bibitem{Jh04} J.\,Johnson. A notion of regularity for differential
local algebras. Contributions to algebra (collection of papers
dedicated to Ellis Kolchin.), pp.~211-232. Academic press, New York,
1977.

\bibitem{Jh02} J.\,Johnson. Differential dimesion polynomials and a fundamental theorem on differential
modules. Amer. J.~Math., 91, 1969, 239-248.

\bibitem{Jh03} J.\,Johnson. Kahler differentials and differential
algebra. Ann. of~Math., (2), 89, 1969, 92-98.

\bibitem{Kl} E.\,R.~Kolchin. Differential Algebra and Algebraic Groups. Academic
Press, New York, 1976.

\bibitem{Mu}H.\,Matsumura. Commutative algebra. The
benjamin/cummings publishing company. 1980.

\bibitem{SW} M.\,van\,der~Put and M.\,F.~Singer. Galois Theory of Linear
Differential Equations  Grundlehren der mathematischen
Wissenschaften, Volume 328, Springer,  2003.

\end{thebibliography}
\end{document}